\documentclass{amsart}
\usepackage{delimset, amssymb, enumitem, mathtools, thmtools}
\usepackage{enumitem}
\setlist{itemsep=4pt, topsep=0pt, leftmargin=17pt}

\usepackage{xcolor}
\definecolor{PKU}{cmyk}{0, 1, 1, .45}
\usepackage[colorlinks, citecolor=blue, linkcolor=PKU, pagebackref]{hyperref}

\usepackage[capitalize]{cleveref}
\crefname{def}{Def.}{Defs.}
\Crefname{def}{Definition}{Definitions}
\creflabelformat{def}{\upshape(#2#1#3)}
\crefname{ineq}{Ineq.}{Ineqs.}
\Crefname{ineq}{Inequality}{Inequalities}
\creflabelformat{ineq}{\upshape(#2#1#3)}

\makeatletter
\newcommand\creflabel[2][\@currentcounter]{%
 \crefalias{\@currentcounter}{#1}\label{#2}}
\makeatother 

\newtheorem{theorem}{Theorem}[section]
\newtheorem{conjecture}[theorem]{Conjecture}
\newtheorem{lemma}[theorem]{Lemma}
\newtheorem{corollary}[theorem]{Corollary}
\newtheorem{proposition}[theorem]{Proposition}
\numberwithin{equation}{section}

\usepackage{tikz}
\usetikzlibrary{arrows.meta}
\tikzset{
edge/.style={semithick},
ball/.style={shape=circle, minimum size=1mm, ball color=black, inner sep=0.5},
ellipsis/.style={shape=circle, fill, inner sep=.5},
range/.style={draw=blue, {Stealth}-{Stealth}}}

\allowbreak
\allowdisplaybreaks

\usepackage[numbers, sort&compress, nonamebreak, merge, elide, longnamesfirst]{natbib}

\author[L. Chen]{L. Chen}
\address[L. Chen]{Department of Mathematics, Southern University of Science and Technology, Shenzhen, Guangdong 518055, China}
\email{chenl2022@mail.sustech.edu.cn}

\author[Y. T. He]{Y. T. He}
\address[Y. T. He]{School of Mathematics and Statistics, Xi'an Jiaotong University, Xi'an, Shaanxi 710049, China}
\email{yitong.he@stu.xjtu.edu.cn}

\author[D. G. L. Wang]{David G. L. Wang}
\address[David G. L. Wang]{School of Mathematics and Statistics \& MIIT Key Laboratory of Mathematical Theory and Computation in Information Security, Beijing Institute of Technology, Beijing 102400, China}
\email{glw@bit.edu.cn}

\title{Clocks are $e$-positive}

\begin{document}

\begin{abstract}
Along with his confirmation of the $e$-positivity of all cycle-chord graphs $\theta_{ab1}$, the third author conjectured the $e$-positivity of all theta graphs $\theta_{abc}$. In this paper, we establish the $e$-positivity of all clock graphs $\theta_{ab2}$ by using the composition method. The key idea is to investigate the fibers of certain partial reversal transformation on compositions with all parts at least $2$.
\end{abstract}
\maketitle
\tableofcontents

\section{Introduction}

In 1995, \citet{Sta95} introduced the concept of 
chromatic symmetric function of a graph $G$ as
\[
X_G(x_1,x_2,\dots)
=\sum_{\kappa\colon V(G)\to\{1,2,\dots\}}
x_{\kappa(v_1)}x_{\kappa(v_2)}\dotsm,
\]
where the sum runs over all proper colorings $\kappa$ of $G$.
It is a generalization of Birkhoff's chromatic
polynomials in the study of the $4$-color problem.
This concept attracted a considerable number of studies 
from algebraic combinatorialists, graph theorists and representation experts.

A leading conjecture in this field is Stanley and Stembridge's conjecture~\cite{SS93}.
In virtue of \citeauthor{Gua13X}'s reduction \cite{Gua13X},
the conjecture can now be stated 
as ``all unit interval graphs are $e$-positive.''
Here a graph $G$ is said to be \emph{$e$-positive} if the expansion of $X_G$ in the basis of elementary symmetric functions has no negative coefficients. 
A close relation between chromatic quasisymmetric functions 
and representations of symmetric groups
on cohomology of regular semisimple Hessenberg varieties
was revealed by \citet{SW16},
see also \citet{BC18} and \citet{Gua16X}.

A natural generalization of Stanley--Stembridge's conjecture
is to characterize all $e$-positive graphs.
The most noted $e$-positive families are complete graphs, paths, cycles,
and graphs with independence number $2$;
they were known early to \citet{Sta95}.
A renowned $e$-positive family is $K$-chains, 
named and proved by \citet{GS01} via a momentous $(e)$-positivity approach.
These are the graphs obtained from a list $(K_{a_1},u_1,v_1)$,
$\dots$, 
$(K_{a_l},u_l,v_l)$ of complete graphs together with two distinct vertices $u_i,v_i\in V(K_{a_i})$,
by identifying $v_i$ and $u_{i+1}$ for all $1\le i\le l-1$.
$K$-chains include, as examples, lollipops and generalized bulls.
More $e$-positive graphs can be found from \cite{AWv24,BCCCGKKLLS24X,CH19,Dah19,Dv18,FHM19,HHT19,LY21,Tom24,Tsu18,TW24X,WW23-JAC}.

There are also research interests on confirming the non-$e$-positivity
of graphs. 
A typical study along this line is \citet{DSv20}'s conjecture
that no tree with maximum degree~$\Delta$ at least $4$ is $e$-positive.
\citet{Zhe22} confirmed it for trees with $\Delta\ge 6$ 
by using elementary and smart combinatorial inequalities
based on \citet{Wol97D}'s connected partition criterion.
Very recently \citet{Tom24X} took a giant step forward 
of solving the $\Delta=5$ case 
by making an intensive combinatorial investigation
on the existence of certain connected partitions.
His method also works for spiders with $\Delta=4$.
In contrast, no rule was known yet to 
determine whether a general spider of $3$ legs is $e$-positive, 
see \citet{WW23-DAM}.

Like cycles and trees,
some graphs have simple structures, radical positions 
and wide applications in graph theory, 
though they are not unit interval in general.
\citet{WZ25} developed a composition method
and confirmed the $e$-positivity of 
hat graphs as an application, which include cycles and tadpoles;
the method was also proved valid for confirming the $e$-positivity of
spiders of the forms $S(4m+2,\,2m+1,\,1)$,
see \citet{TWW24X}.
It was a conjecture of Aliniaeifard, Wang and van Willigenburg,
see \citet[Conjecture 6.2]{Zhe22}.

To confirm the $e$-positivity of a graph family often requires specific algebraic combinatorial skills. Besides the above,
common skills include the generating function approach,
constructing recurrence relations and sign-reversing involutions.
An incomplete list of proof strategies can be found from \citet{QTW24X}.

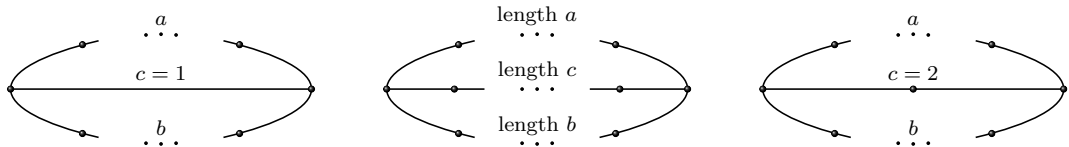
\begin{figure}[htbp]
\begin{tikzpicture}
\coordinate (left) at (180: 2);
\coordinate (left-u) at (180-29.5: 1.2);
\coordinate (left-d) at (180+29.5: 1.2);
\coordinate (c) at (0, 0);
\coordinate (right) at (0: 2);
\coordinate (right-u) at (29.5: 1.2);
\coordinate (right-d) at (-29.5: 1.2);

\coordinate (up) at (0, .73);
\coordinate (up-l) at (-.2, .72);
\coordinate (up-r) at (.2, .72);
\coordinate (down) at (0, -.73);
\coordinate (down-l) at (-.2, -.72);
\coordinate (down-r) at (.2, -.72);

\draw[edge] (left) -- (right);
\draw[edge] (left) arc [start angle=180, delta angle=-40, x radius=5, y radius=1];
\draw[edge] (left) arc [start angle=180, delta angle=40, x radius=5, y radius=1];
\draw[edge] (right) arc [start angle=0, delta angle=-40, x radius=5, y radius=1];
\draw[edge] (right) arc [start angle=0, delta angle=40, x radius=5, y radius=1];

\foreach \b in {left, left-u, left-d, right, right-u, right-d}
\node[ball] at (\b) {};

\foreach \p in {up, up-l, up-r, down, down-l, down-r}
\node[ellipsis] at (\p) {};

\node[above, font=\footnotesize] at (up) {$a$};
\node[above, font=\footnotesize] at (c) {$c=1$};
\node[above, font=\footnotesize] at (-90: .73) {$b$};

\begin{scope}[xshift=5cm]
\coordinate (left) at (180: 2);
\coordinate (left-u) at (180-29.5: 1.2);
\coordinate (left-c) at (180: 1.1);
\coordinate (left-d) at (180+29.5: 1.2);
\coordinate (right) at (0: 2);
\coordinate (right-u) at (29.5: 1.2);
\coordinate (right-c) at (0: 1.1);
\coordinate (right-d) at (-29.5: 1.2);

\coordinate (up) at (0, .73);
\coordinate (up-l) at (-.2, .72);
\coordinate (up-r) at (.2, .72);
\coordinate (c) at (0, 0);
\coordinate (c-l) at (-.2, 0);
\coordinate (c-r) at (.2, 0);
\coordinate (down) at (0, -.73);
\coordinate (down-l) at (-.2, -.72);
\coordinate (down-r) at (.2, -.72);

\draw[edge] (left) -- (-.7, 0);
\draw[edge] (left) arc [start angle=180, delta angle=-40, x radius=5, y radius=1];
\draw[edge] (left) arc [start angle=180, delta angle=40, x radius=5, y radius=1];

\draw[edge] (right) -- (.7, 0);
\draw[edge] (right) arc [start angle=0, delta angle=-40, x radius=5, y radius=1];
\draw[edge] (right) arc [start angle=0, delta angle=40, x radius=5, y radius=1];

\foreach \b in {left, left-u, left-c, left-d, right, right-u, right-c, right-d}
\node[ball] at (\b) {};

\foreach \p in {up, up-l, up-r, c, c-l, c-r, down, down-l, down-r}
\node[ellipsis] at (\p) {};

\node[above, font=\footnotesize] at (up) {length $a$};
\node[above, font=\footnotesize] at (c) {length $c$};
\node[above, font=\footnotesize] at (down) {length $b$};
\end{scope}

\begin{scope}[xshift=10cm]
\coordinate (left) at (180: 2);
\coordinate (left-u) at (180-29.5: 1.2);
\coordinate (left-d) at (180+29.5: 1.2);
\coordinate (c) at (0, 0);
\coordinate (right) at (0: 2);
\coordinate (right-u) at (29.5: 1.2);
\coordinate (right-d) at (-29.5: 1.2);

\coordinate (up) at (0, .73);
\coordinate (up-l) at (-.2, .72);
\coordinate (up-r) at (.2, .72);
\coordinate (down) at (0, -.73);
\coordinate (down-l) at (-.2, -.72);
\coordinate (down-r) at (.2, -.72);

\draw[edge] (left) -- (right);
\draw[edge] (left) arc [start angle=180, delta angle=-40, x radius=5, y radius=1];
\draw[edge] (left) arc [start angle=180, delta angle=40, x radius=5, y radius=1];
\draw[edge] (right) arc [start angle=0, delta angle=-40, x radius=5, y radius=1];
\draw[edge] (right) arc [start angle=0, delta angle=40, x radius=5, y radius=1];

\foreach \b in {left, left-u, left-d, right, right-u, right-d, c}
\node[ball] at (\b) {};

\foreach \p in {up, up-l, up-r, down, down-l, down-r}
\node[ellipsis] at (\p) {};

\node[above, font=\footnotesize] at (up) {$a$};
\node[above, font=\footnotesize] at (c) {$c=2$};
\node[above, font=\footnotesize] at (down) {$b$};
\end{scope}

\end{tikzpicture}
\caption{A cycle-chord, a theta graph, and a clock.}
\label{fig:theta}
\end{figure}

Recently the third author \cite{Wang25} 
established the $e$-positivity of cycle-chords,
which are obtained by identifying an edge of two cycles,
see the left figure in \cref{fig:theta}.
He realized that cycle-chords 
can be regarded as obtained by connecting two distinct vertices 
with three disjoint paths and by restricting one of them to have length $1$.
If one uses $k$ disjoint paths to connect two vertices,
then the resulting graphs are paths if $k=1$, and cycles if $k=2$.
When $k\ge 4$, he found some infinite families of graphs that are not $e$-positive.
For $k=3$, the resulting graphs were named 
\emph{theta graphs} by \citet{Bon72}. To be more precise,
a theta graph is the union of three internally disjoint paths 
that have the same two distinct ends.
When the three paths have lengths $a$, $b$ and $c$,
we denote the theta graph by $\theta_{abc}$.
see the middle figure in \cref{fig:theta}.
The $e$-positivity of cycle-chords then inspires \citeauthor{Wang25} 
to pose the following conjecture.
\begin{conjecture}[\citeauthor{Wang25}]\label{conj:theta}
All theta graphs $\theta_{abc}$ are $e$-positive.
\end{conjecture}

In this paper,
we confirm \cref{conj:theta} for $c=2$.
Imaging the cycle of length $a+b$ in $\theta_{ab2}$ as a clock face, 
and the two edges that constitute the length $2$ path 
as the hour hand and minute hand,
we call the theta graphs $\theta_{ab2}$
\emph{clock graphs}, see the right figure in \cref{fig:theta}.
Here is our main result.

\begin{theorem}\label{thm:X.clock}
Clocks are $e$-positive.
\end{theorem}

This paper is organized as follows.
In \cref{sec:pre},
we give an overview of basic knowledge on chromatic symmetric functions,
the composition method, as well as positive $e_I$-expansions
of some particular graphs that will be of use.
In \cref{sec:theta}, we derive two $e_I$-expansions 
for the chromatic symmetric function of all theta graphs.
In \cref{sec:clock},
we establish \cref{thm:X.clock} based on one of the $e_I$-expansions.
The main idea in our proof is to figure out the 
fibers of certain partial reversal transformation on compositions 
with all parts at least $2$.

\section{Preliminaries}\label{sec:pre}
This section contains 
basic knowledge on chromatic symmetric functions
that will be of use.
We adopt terminologies from \citet{GKLLRT95} and \citet{Sta11B}.
Let $n$ be a positive integer. 
A \emph{composition} of $n$ is 
a sequence of positive integers with sum~$n$,
commonly denoted 
$I=i_1 \dotsm i_z\vDash n$, 
with \emph{size} $\abs{I}=n$, \emph{length} $\ell(I)=z$, 
and \emph{parts} $i_1,\dots,i_z$.
When all parts $i_k$ have the same value $i$, we write $I=i^z$.
The \emph{reversal} composition $i_z\dotsm i_1$ is written as $\overline{I}$.
For convenience, we denote the composition obtained 
by removing the $k$th part by~$I\backslash i_k$, i.e.,
\[
I\backslash i_k
=i_1\dotsm i_{k-1}i_{k+1}\dotsm i_z.
\]
When a capital letter like $I$ or $J$
stands for a composition,
its small letter counterpart with integer subscripts stands for the parts.
A \emph{partition}\index{partition} of $n$
is a multiset of positive integers $\lambda_i$ with sum~$n$,
denoted 
$\lambda=\lambda_1\lambda_2\dotsm
\vdash n$,
where $\lambda_1\ge \lambda_2\ge\dotsm\ge 1$.

A homogeneous \emph{symmetric function} of degree $n$ over
the field $\mathbb{Q}$ of rational numbers is a formal power series
\[
f(x_1, x_2, \dots)
=\sum_{
\lambda}
c_\lambda \cdotp
x_1^{\lambda_1} 
x_2^{\lambda_2} 
\dotsm,
\]
where the sum is over all weak compositions $\lambda$
(with $0$ as parts allowed),
such that $f(x_1, x_2, \dots)=f(x_{\pi(1)}, x_{\pi(2)}, \dots)$
for any permutation $\pi$.
Let $\operatorname{Sym}^0=\mathbb{Q}$, and let $\operatorname{Sym}^n$ 
be the vector space of homogeneous symmetric functions of degree $n$ over~$\mathbb{Q}$. 
One basis of $\operatorname{Sym}^n$ consists of elementary symmetric functions 
$e_\lambda$ for all partitions $\lambda\vdash n$, where
\[
e_\lambda
=e_{\lambda_1}
e_{\lambda_2}\dotsm
\quad\text{and}\quad
e_k
=\sum_{1\le i_1<\dots<i_k} 
x_{i_1} \dotsm x_{i_k}.
\]
A symmetric function $f\in\mathrm{Sym}$ is said to be \emph{$e$-positive}
if every $e_\lambda$-coefficient of $f$ is nonnegative.

\Citet{Sta95} introduced the chromatic symmetric function 
for a simple graph $G=(V,E)$ as
\[
X_G
=
\sum_{\kappa\colon V\to\{1,2,\dots\}}
\prod_{v\in V}
x_{\kappa(v)},
\]
where $\kappa$ runs over proper colorings of~$G$. 
One of the most popular tools in studying chromatic symmetric functions is 
the triple-deletion property established by
\citet[Theorem~3.1, Corollaries 3.2 and 3.3]{OS14}.

\begin{proposition}[\citeauthor{OS14}]\label{prop:3del}
Let $G$ be a graph with a stable set $T$ of cardinality $3$.
Denote by $e_1'$, $e_2'$ and $e_3'$
the edges linking the vertices in $T$.
For any set $S\subseteq \{1,2,3\}$, 
denote by $G_S$ the graph with vertex set~$V(G)$
and edge set $E(G)\cup\{e_j'\colon j\in S\}$.
Then 
\[
X_{G_{12}}
=X_{G_1}+X_{G_{23}}-X_{G_3}
\quad\text{and}\quad
X_{G_{123}}
=X_{G_{13}}+X_{G_{23}}-X_{G_3}.
\]
\end{proposition}

For any composition $I$, there is a unique partition $\rho(I)$
which consists of the parts of $I$.
This allows us to define 
$e_I=e_{\rho(I)}$.
An \emph{$e_I$-expansion} of a symmetric function $f\in\mathrm{Sym}^n$ is an expression
$f=\sum_{I\vDash n} c_I e_I$.
We call it a \emph{positive $e_I$-expansion}
if $c_I\ge 0$ for all $I$.
\citet[Table~1]{SW16} discovered 
a captivating positive $e_I$-expansion for paths $P_n$.

\begin{proposition}[\citeauthor{SW16}]\label{prop:X.path}
We have $X_{P_n}
=\sum_{I\vDash n}
w_I
e_I$ for $n\ge 1$, where
\begin{equation}\label{def:w}
w_I
=i_1(i_2-1)(i_3-1)\dotsm(i_z-1)
\quad\text{if $I=i_1i_2\dotsm i_z$}.
\end{equation}
\end{proposition}
Analogously, \citet[Corollary 6.2]{Ell17X} 
gave a formula for the chromatic quasisymmetric function of cycles $C_n$, whose
$t=1$ specialization is the following.

\begin{proposition}[\citeauthor{Ell17X}]\label{prop:X.cycle}
We have 
$
X_{C_n}
=\sum_{I\vDash n}
(i_1-1)
w_I
e_I$ for $n\ge 2$.
\end{proposition}

For the purpose of establishing the $e$-positivity of a graph $G$,
it suffices to 
present a positive $e_I$-expansion for the chromatic symmetric function of $G$.
This approach was developed by \citet{WZ25},
called the \emph{composition method}.
Its first application was for tadpole graphs.
The \emph{tadpole}~$C_a^l$ is the graph obtained by identifying
a vertex on the cycle $C_a$ and an end of the path $P_{l+1}$. 
It has the same number $a+l$ of vertices and edges.
\citet{WZ25} expressed the chromatic symmetric function of tadpoles 
in terms of the function
\[
\Theta_{I}^{+}(a)=\sigma_{I}^{+}(a)-a,
\quad
\text{
where
$
\sigma_{I}^{+}(a)
=\min\{
i_1+\dots+i_k
\colon
0\le k\le \ell(I),\
i_1+\dots+i_k\ge a\}$.
}
\]

\begin{theorem}[\citeauthor{WZ25}]\label{thm:X.tadpole}
We have 
$
X_{C_{n-l}^l}
=\sum_{I\vDash n}
\Theta_I^+(l+1)
w_I
e_I$
for any $0\le l\le n-2$.
\end{theorem}

We remark that \cref{thm:X.tadpole} reduces to \cref{prop:X.path}
for when $l=n-2$, and to \cref{prop:X.cycle} for when $l=0$.
\citet{WZ25} also defined the function
\[
\Theta_I^-(a)=a-\sigma_{I}^-(a),
\quad\text{where $
\sigma_I^-(a)
=\max\{
i_1+\dots+i_k
\colon
0\le k\le \ell(I),\
i_1+\dots+i_k\le a
\}$},
\]
and presented the relationship between
the functions $\Theta_I^-$ and $\Theta_I^+$ as follows,
see \cite[Lemma 2.9]{WZ25}.

\begin{lemma}[\citeauthor{WZ25}]\label{lem:Theta+-}
For any composition $I\vDash n$ and any real number $0\le a\le n$, 
\[
\Theta_{I}^-(a)=\Theta_{\overline{I}}^{+}(n-a).
\]
\end{lemma}

Cycle-chord graphs
$\mathrm{CC}_{a,b}$
are the theta graphs $\theta_{ab1}$, see \cref{fig:theta}. 
\citet{Wang25} provided a positive $e_I$-expansion for cycle-chords,
as well as an intuitive interpretation of the $e_I$-coefficients.

\begin{theorem}[\citeauthor{Wang25}]\label{thm:X.CC}
For any $a,b\geq 2$ and $n=a+b$,  
\[
X_{\mathrm{CC}_{a,b}}
=\sum_{I\vDash n}
\brk4{
\sum_{i=1}^b\Theta_{I}^+(i)
-\sum_{i=1}^{b-1}\Theta_{\overline{I}}^-(i)}
w_I
e_I
=\sum_{I=i_1i_2\dotsm i_z\vDash n}
\Delta_{I}(b)w_{I}e_{I},
\]
where
\begin{equation}\creflabel[def]{def:Delta}
\Delta_I(b)
=\begin{dcases*}
s(i_p-s-i_1),
& if $i_1\le i_p-s$,\\
e_2(i_p-s,\,i_{p+1},\,\dots,\,i_q,\,t),
& otherwise,
\end{dcases*}
\end{equation}
in which $e_2(x_1,\dots,x_m)=\sum_{1\le i<j\le m}x_i x_j$,
and the symbols $p,q,s,t$ are defined by 
\[
b=i_1+\dots+i_{p-1}+s
=i_2+\dots+i_q+t,
\]
such that $1\le p,q\le z$, $1\le s\le i_p$, and $1\le t\le i_{q+1}$,
with the convention $i_{z+1}=i_1$.
As a consequence, all cycle-chords are $e$-positive.
\end{theorem}

Here we copy from \citet[Lemma 3.2]{Wang25}
some relations among the numbers $p$, $q$, $s$ and $t$.
They will be of frequent use
in \cref{sec:theta,sec:clock}.

%: lem:psqt
\begin{lemma}[\citeauthor{Wang25}]\label{lem:psqt}
Let $I\vDash n$. Then we have the following.
\begin{enumerate}
\item\label[itm]{itm:ip-s}
$i_p-s=\Theta_{\overline{I}}^-(a)$.
\item\label[itm]{itm:q>=p-1}
$q\ge p-1$, and
$i_p+i_{p+1}+\dots+i_q
=i_1+s-t$.
\item\label[itm]{itm:q=p-1:iff}
$q=p-1\iff i_1\le i_p-s$.
\item\label[itm]{itm:q=z:iff}
$a-i_1
=
i_{q+1}+\dots+i_z-t$.
As a consequence,
$q=z\iff i_1>a$.
\end{enumerate}
\end{lemma}

\section{Two $e_I$-expansions for theta graphs}\label{sec:theta}
We start from manufacturing an $e_I$-expansion for $X_{\theta_{abc}}$.

\begin{theorem}\label{thm:X.theta.abc}
Let $a\ge b\ge c\ge 1$ and $n=a+b+c-1$.
Then 
\[
X_{\theta_{abc}}
=
\sum_{I\vDash n}
\brk1{
d_I+\Delta_I(b+c-1)
}
w_I
e_I,
\]
where
\[
d_I
=
\sum_{k=2}^{c}\Theta_{I}^{+}(k)
-\sum_{k=a}^{a+c-2}\Theta_{\overline{I}}^-(k).
\]
\end{theorem}
\begin{proof}
Let $G=\theta_{abc}$. 
Applying \cref{prop:3del} for the edges~$e_i$ illustrated in \cref{fig:theta-tri}, we obtain
\[
X_{\theta_{abc}}
=X_{\theta_{(a+1)b(c-1)}}
+X_{C_{a+b}^{c-1}}
-X_{C_{c+b-1}^{a}},
\quad\text{for $c\ge 2$.}
\]
\begin{figure}[htbp]
\begin{tikzpicture}
\coordinate (left) at (180: 2);
\coordinate (left-u) at (180-29.5: 1.2);
\coordinate (left-c) at (180: 1.1);
\coordinate (left-d) at (180+29.5: 1.2);
\coordinate (right) at (0: 2);
\coordinate (right-u) at (29.5: 1.2);
\coordinate (right-c) at (0: 1.1);
\coordinate (right-d) at (-29.5: 1.2);

\coordinate (up) at (0, .73);
\coordinate (up-l) at (-.2, .72);
\coordinate (up-r) at (.2, .72);
\coordinate (c) at (0, 0);
\coordinate (c-l) at (-.2, 0);
\coordinate (c-r) at (.2, 0);
\coordinate (down) at (0, -.73);
\coordinate (down-l) at (-.2, -.72);
\coordinate (down-r) at (.2, -.72);

\draw[edge] (left) -- (-.7, 0);
\draw[edge] (left) arc [start angle=180, delta angle=-40, x radius=5, y radius=1];
\draw[edge] (left) arc [start angle=180, delta angle=40, x radius=5, y radius=1];

\draw[edge] (right) -- (.7, 0);
\draw[edge] (right) arc [start angle=0, delta angle=-40, x radius=5, y radius=1];
\draw[edge] (right) arc [start angle=0, delta angle=40, x radius=5, y radius=1];

\draw[PKU, thick] 
(left-c) 
-- (left) 
node[midway, above, font=\footnotesize]{$e_2$}
arc[start angle=180, delta angle=36, x radius=5, y radius=1] 
node[midway, below, font=\footnotesize]{$e_1$};
\draw[PKU, thick, densely dotted] 
(left-c) -- (left-d) node[midway, right, font=\footnotesize]{$e_3$};

\foreach \b in {left-u, right, right-u, right-c, right-d}
\node[ball] at (\b) {};

\foreach \b in {left, left-c, left-d}
\node[ball, ball color=PKU] at (\b) {};

\foreach \p in {up, up-l, up-r, c, c-l, c-r, down, down-l, down-r}
\node[ellipsis] at (\p) {};

\node[above, font=\footnotesize] at (up) {$a$};
\node[above, font=\footnotesize] at (c) {$b$};
\node[above, font=\footnotesize] at (down) {$c$};
\end{tikzpicture}
\caption{Triple deletion of the theta graph $\theta_{abc}$.}
\label{fig:theta-tri}
\end{figure}
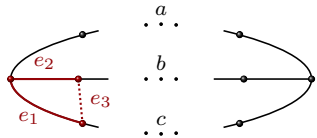

Using it iteratively, we can deduce that
\[
X_G
=
X_{\theta_{(n-b)b1}}
+\sum_{k=1}^{c-1}X_{C_{a+b+k-1}^{c-k}}
-\sum_{k=1}^{c-1}X_{C_{b+c-k}^{a+k-1}}.
\]
Since $\theta_{(n-b)b1}=\mathrm{CC}_{(n-b)b}$,
we can deduce using \cref{thm:X.CC,thm:X.tadpole} that
\begin{equation}\label{pf:X.theta.abc}
X_G
=\sum_{I\vDash n}
\brk4{
\sum_{k=1}^{b}
\Theta_{I}^+(k)
-
\sum_{k=1}^{b-1}
\Theta_{\overline{I}}^-(k)}
w_I
e_I
+
\sum_{k=1}^{c-1}
\sum_{I\vDash n}
\brk4{
\Theta_{I}^{+}(c-k+1)
-
\Theta_{I}^{+}(a+k)
}
w_I
e_I.
\end{equation}
We regard \cref{pf:X.theta.abc} as
the sum of four double sums.
The second last double sum can be rewritten as
\begin{equation}\label{pf:X.theta.abc.3}
\sum_{k=1}^{c-1}
\sum_{I\vDash n}
\Theta_{I}^{+}(c-k+1)
=
\sum_{I\vDash n}
\sum_{k=2}^{c}
\Theta_{I}^{+}(k).
\end{equation}
The last double sum in \cref{pf:X.theta.abc}
can be rewritten by using \cref{lem:Theta+-} as
\[
-
\sum_{k=1}^{c-1}
\sum_{I\vDash n}
\Theta_{I}^{+}(a+k)
w_I
e_I
=
-
\sum_{I\vDash n}
\sum_{k=1}^{c-1}
\Theta_{\overline{I}}^-(n-a-k)
w_I
e_I
=
-
\sum_{I\vDash n}
\sum_{k=b}^{b+c-2}
\Theta_{\overline{I}}^-(k)
w_I
e_I,
\]
which can then be merged with the second double sum in \cref{pf:X.theta.abc} as 
\begin{equation}\label{pf:X.theta.abc.24}
-\sum_{I\vDash n}
\sum_{k=1}^{b-1}
\Theta_{\overline{I}}^-(k)
w_I
e_I
-\sum_{k=1}^{c-1}
\sum_{I\vDash n}
\Theta_{I}^{+}(a+k)
w_I
e_I
=
-\sum_{I\vDash n}
\sum_{k=1}^{b+c-2}
\Theta_{\overline{I}}^-(k)
w_I
e_I.
\end{equation}
Substituting \cref{pf:X.theta.abc.3,pf:X.theta.abc.24} back into \cref{pf:X.theta.abc}, we obtain a simplification
\begin{equation}\label{pf:X.theta.abc.new}
X_G
=\sum_{I\vDash n}
\brk4{
\sum_{k=1}^{b}
\Theta_{I}^+(k)
+
\sum_{k=2}^{c}
\Theta_{I}^+(k)
-\sum_{k=1}^{b+c-2}
\Theta_{\overline{I}}^-(k)
}
w_I
e_I.
\end{equation}
In view of \cref{thm:X.CC,thm:X.CC}, 
the negative term in \cref{pf:X.theta.abc.new} appears in
\[
\sum_{I\vDash n}
\brk4{
\sum_{k=1}^{b+c-1}
\Theta_{I}^+(k)
-
\sum_{k=1}^{b+c-2}
\Theta_{\overline{I}}^-(k)
}
w_I
e_I
=X_{\mathrm{CC}_{a(b+c-1)}}
=\sum_{I\vDash n}
\Delta_I(b+c-1)
w_I
e_I.
\] 
We therefore separate the first positive term in \cref{pf:X.theta.abc.new} as
\[
\sum_{k=1}^b\Theta_I^+(k)
=\sum_{k=1}^{b+c-1}\Theta_I^+(k)
-\sum_{k=b+1}^{b+c-1}\Theta_I^+(k),
\]
and transform its last term by using \cref{lem:Theta+-} again as 
\[
-
\sum_{k=b+1}^{b+c-1}
\Theta_{I}^+(k)
w_I
e_I
=
-
\sum_{k=b+1}^{b+c-1}
\Theta_{\overline{I}}^-(n-k)
w_I
e_I
=
-
\sum_{k=a}^{a+c-2}
\Theta_{\overline{I}}^-(k)
w_I
e_I.
\]
The desired formula follows by
integrating these relations.
\end{proof}

%: def:varphi
Let $I=PQ$ be a composition, where $Q$ 
is the shortest suffix of $I$ 
with size at least $a$. Define
\begin{equation}\creflabel[def]{def:phi}
\varphi(I)
=\begin{dcases*}
I,
& if $Q=I$, \\
i_1\overline{P\backslash i_1}Q,
& if $Q\ne I$.
\end{dcases*}
\end{equation}
It is clear that 
$w_Ie_I=w_{\varphi(I)}e_{\varphi(I)}$ and 
that $\varphi$ is an involution, i.e., $
\varphi(\varphi(I))=I$.
For any composition~$I$ of $a+b+c-1$, define
\begin{equation}\creflabel[def]{def:d'}
d_I'
=
\sum_{k=2}^{c}
\Theta_{I}^{+}(k)
-\sum_{k=a}^{a+c-2}
\Theta_{\overline{\varphi(I)}}^-(k).
\end{equation}

\begin{theorem}\label{thm:theta.abc:phi}
Let $n=a+b+c-1$, where $a\ge b\ge c\ge 1$.
Then 
\[
X_{\theta_{abc}}
=\sum_{I\vDash n}
\brk1{d_I'+\Delta_I(b+c-1)
}
w_I
e_I.
\]
\end{theorem}
\begin{proof}
Recall the definition of $d_I$ from \cref{thm:X.theta.abc}.
Since $\varphi$ is an involution, 
\[
d_I+d_{\varphi(I)}
=
\sum_{k=2}^{c}\Theta_{I}^{+}(k)
-\sum_{k=a}^{a+c-2}\Theta_{\overline{I}}^-(k)
+
\sum_{k=2}^{c}\Theta_{\varphi(I)}^{+}(k)
-\sum_{k=a}^{a+c-2}\Theta_{\overline{\varphi(I)}}^-(k)
=
d_I'+d_{\varphi(I)}'.
\]
It follows that 
\[
\sum_{I\vDash n}
d_I
w_Ie_I
=\frac{1}{2}
\sum_{I\vDash n}
\brk1{
d_I
+d_{\varphi(I)}
}
w_Ie_I
=\frac{1}{2}
\sum_{I\vDash n}
\brk1{
d_I'
+d_{\varphi(I)}'
}
w_Ie_I
=
\sum_{I\vDash n}
d_I'
w_Ie_I.
\]
By \cref{thm:X.theta.abc}, we obtain the desired formula.
\end{proof}

\begin{corollary}\label{cor:theta2:DI}
Let $a\ge b\ge 2$. Then $X_{\theta_{ab2}}
=\sum_{I\vDash n}D_Iw_Ie_I$,
where
\begin{equation}\creflabel[def]{def:DI}
D_I
=\Theta_I^+(2)
-\Theta_{\overline{\varphi(I)}}^-(a)
+\Delta_I(b+1).
\end{equation}
\end{corollary}
\begin{proof}
This is immediate by taking $c=2$ in \cref{thm:theta.abc:phi}.
\end{proof}

Before ending this section, 
we cope with the compositions $I$ with a suffix of size $a$.
Let
\[
\mathcal A
=\brk[c]1{
I\vDash n
\colon
w_I>0,\
\Theta_{\overline{I}}^{+}(a)=0
}.
\]
Then the map restriction $\varphi|_{\mathcal A}$ is an involution on $\mathcal A$.
Along the same lines, one may infer that 
\[
\sum_{I\in \mathcal A}
\brk1{d_I+\Delta_I(b+c-1)}
w_I
e_I
=
\sum_{I\in \mathcal A}
\brk1{d_I'+\Delta_I(b+c-1)}
w_I
e_I.
\]

\begin{lemma}\label{lem:A}
Let $a\ge b\ge c\ge 2$ and $I\in\mathcal A$.
Then $d_I'\ge 0$,
where $d_I'$ is defined by \cref{def:d'}.
\end{lemma}
\begin{proof}
Let $I=i_1\dotsm i_z\in \mathcal A$.
Then there exists a unique part index $1\le p\le z-1$ such that 
\begin{equation}\creflabel[def]{def:p:suffix=a}
a=
i_{p+1}+\dots+i_z.
\end{equation}
Since $b\ge c\ge 2$, we find
\[
i_1+\dots+i_p
=b+c-1
\ge c+1.
\]
In order to simplify the two sums in \cref{def:d'} of $d_I'$, 
we need to identify some special parts of $I$.
Recall that 
\[
\sigma_I^-(a)
=\max\{
i_1+\dots+i_k
\colon
0\le k\le \ell(I),\
i_1+\dots+i_k\le a
\}.
\]
Define $m$ and $t$ to be the indices such that 
\AddToHook{env/align/begin}{\crefalias{equation}{def}}
\begin{align}
\label{def:m}
\sigma_I^-(c-1)
&=i_1+\dots+i_m
\quad\text{and}\\
\label{def:t}
\sigma_{i_2i_3\dotsm i_p}^-(c-2)
&=i_2+\dots+i_t.
\end{align}
Then $0\le m\le p-1$ and $1\le t\le p$. 
In particular, 
\begin{itemize}
\item
$m=0\iff i_1\ge c$, and
\item
$t=1\iff$ (i) $p=1$ and $i_1=b+c-1$,  or (ii) $p\ge 2$ and $i_2\ge c-1$.
\end{itemize}
We claim that 
$m\le t$.
In fact, if $m>t$,
then $m\ge t+1\ge 2$.
By \cref{def:m}, we deduce that 
\[
c-1
\ge
i_1+\dots+i_m
\ge
i_1+(i_2+\dots+i_{t+1})
>
i_1+(c-2)
\ge c-1,
\]
a contradiction. 
This proves the claim.

Next, we define
\AddToHook{env/align/begin}{\crefalias{equation}{def}}
\begin{align}
\label{pf:d>=0:def:u}
u&=c-(i_1+\dots+i_m)
\quad\text{and}\\
\label{pf:d>=0:def:v}
v&=c-2-(i_2+\dots+i_t).
\end{align}
We now seek for bounds of $u$ and $v$.
By \cref{def:m}, 
we find
\begin{equation}\creflabel[ineq]{ub:u}
1\le
u\le 
(i_1+\dots+i_{m+1})
-(i_1+\dots+i_m)
=i_{m+1}.
\end{equation}
On the other hand, by \cref{def:t}, we find $v\ge 0$. 
For an upper bound of $v$, we argue as follows.
\begin{enumerate}
\item
If $t\le p-1$, then $p\ge t+1\ge 2$.
In this case, \cref{def:t} implies
\begin{equation}\creflabel[ineq]{ub:v:t<=p-1}
v
\le
i_{t+1}-1.
\end{equation}
\item
If $t=p$, 
then by the premise $i_1+\dots+i_p=b+c-1$ we deduce that 
\[
v
=
i_1+\dots+i_p-b+1
-2
-(i_2+\dots+i_t)
=i_1-b-1
\le i_1-3.
\]
\end{enumerate}
See \cref{fig:I.a} for the positions of the parameters $p$, $m$, $t$, $u$ and $v$ in $I$, 
and \cref{fig:reverse.phiI.a} for those in $\overline{\varphi(I)}$.
\begin{figure}[htbp]
\begin{tikzpicture}
\draw (0, 0) -- (15, 0);
\draw (0, 0) -- (0, .2);
\draw (1, 0) -- (1, .2);
\draw (2, 0) -- (2, .2);
\draw (3, 0) -- (3, .2);
\draw (4, 0) -- (4, .2);
\draw (6, 0) -- (6, .2);
\draw (7, 0) -- (7, .2);
\draw (8, 0) -- (8, .2);
\draw (10, 0) -- (10, .2);
\draw (11, 0) -- (11, .2);
\draw (12, 0) -- (12, .2);
\draw (13, 0) -- (13, .2);
\draw (14, 0) -- (14, .2);
\draw (15, 0) -- (15, .2);

\node[above] at (0.5, 0){$i_1$};
\node[above] at (1.5, 0){$i_2$};
\node[above] at (2.5, 0){$\dotsm$};
\node[above] at (3.5, 0){$i_m$};
\node[above] at (5, 0){$i_{m+1}$};
\node[above] at (6.5, 0){$\dotsm$};
\node[above] at (7.5, 0){$i_t$};
\node[above] at (9, 0){$i_{t+1}$};
\node[above] at (10.5, 0){$\dotsm$};
\node[above] at (11.5, 0){$i_p$};
\node[above] at (12.5, 0){$i_{p+1}$};
\node[above] at (13.5, 0){$\dotsm$};
\node[above] at (14.5, 0){$i_z$};

\draw[draw=blue] (12, -1) -- (12, 0);
\draw[draw=blue] (15, -1) -- (15, 0);
\draw[draw=blue, {Stealth}-{Stealth}]
(12, -.6) -- (13.5, -.6) node[fill=white]{$a$}--(15, -.6);

\draw[draw=blue] (0, -1) -- (0, 0);
\draw[draw=blue] (5, -1) -- (5, 0);
\draw[draw=blue, {Stealth}-{Stealth}]
(0, -.6) -- (2.5, -.6) node[fill=white]{$c$}--(5, -.6);

\draw[draw=blue] (1, 0) -- (1,1);
\draw[draw=blue] (9, 0) -- (9,1);
\draw[draw=blue, {Stealth}-{Stealth}]
(1, .7) -- (5, .7) node[fill=white]{$c-2$}--(9, .7);

\draw[draw=red] (8, 0) -- (8, -.6);
\draw[draw=red] (9, 0) -- (9, -.6);
\draw[draw=red, {Stealth}-{Stealth}]
(8, -.3) -- (8.5, -.3) node[fill=white]{$v$}--(9, -.3);

\draw[draw=blue] (5, 0) -- (5, -1);
\draw[draw=blue] (12, 0) -- (12, -1);
\draw[draw=blue, {Stealth}-{Stealth}]
(5, -.6) -- (7, -.6) node[fill=white]{$b-1$}--(12, -.6);
\draw[draw=red] (5, -.6) -- (5, 0);
\draw[draw=red] (4, -.6) -- (4, 0);
\draw[draw=red,{Stealth}-{Stealth}]
(4, -.3) -- (4.5, -.3) node[fill=white]{$u$}--(5, -.3);
\end{tikzpicture}
\caption{The numbers $p$, $m$, $t$, $u$ and $v$ for $I\in\mathcal A$ that are defined by \cref{def:p:suffix=a,def:m,def:t,pf:d>=0:def:u,pf:d>=0:def:v} respectively.}
\label{fig:I.a}
\end{figure}
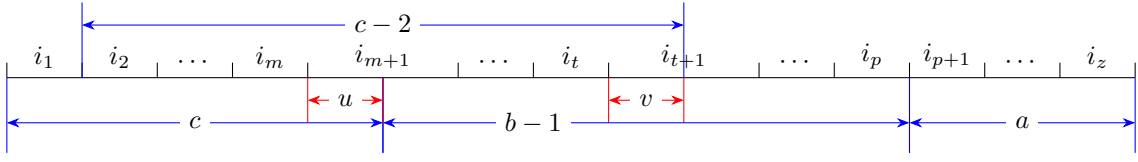

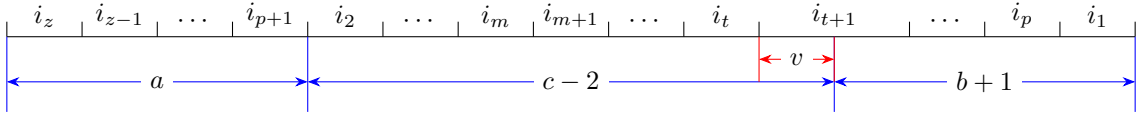
\begin{figure}[htbp]
\begin{tikzpicture}
\draw (0, 0) -- (15, 0);
\draw (0, 0) -- (0, .2);
\draw (1, 0) -- (1, .2);
\draw (2, 0) -- (2, .2);
\draw (3, 0) -- (3, .2);
\draw (4, 0) -- (4, .2);
\draw (5, 0) -- (5, .2);
\draw (6, 0) -- (6, .2);
\draw (7, 0) -- (7, .2);
\draw (8, 0) -- (8, .2);
\draw (9, 0) -- (9, .2);
\draw (10, 0) -- (10, .2);

\draw (12, 0) -- (12, .2);
\draw (13, 0) -- (13, .2);
\draw (14, 0) -- (14, .2);
\draw (15, 0) -- (15, .2);

\node[above] at (0.5, 0){$i_z$};
\node[above] at (1.5, 0){$i_{z-1}$};
\node[above] at (2.5, 0){$\dotsm$};

\node[above] at (3.5, 0){$i_{p+1}$};
\node[above] at (4.5, 0){$i_{2}$};
\node[above] at (5.5, 0){$\dotsm$};
\node[above] at (6.5, 0){$i_m$};
\node[above] at (7.5, 0){$i_{m+1}$};
\node[above] at (8.5, 0){$\dotsm$};
\node[above] at (9.5, 0){$i_t$};

\node[above] at (11, 0){$i_{t+1}$};
\node[above] at (12.5, 0){$\dotsm$};
\node[above] at (13.5, 0){$i_p$};
\node[above] at (14.5, 0){$i_1$};

\draw[draw=blue] (0, 0) -- (0, -1);
\draw[draw=blue] (4, -1) -- (4, 0);
\draw[draw=blue, {Stealth}-{Stealth}]
(0, -.6) -- (2, -.6) node[fill=white]{$a$}--(4, -.6);

\draw[draw=blue] (11, -1) -- (11, 0);
\draw[draw=blue, {Stealth}-{Stealth}]
(4, -.6) -- (7.5, -.6) node[fill=white]{$c-2$}--(11, -.6);

\draw[draw=red] (10, -.6) -- (10, 0);
\draw[draw=red] (11, -.6) -- (11, 0);
\draw[draw=red,{Stealth}-{Stealth}]
(10, -.3) -- (10.5, -.3) node[fill=white]{$v$} -- (11, -.3);

\draw[draw=blue] (15, -1) -- (15, 0);
\draw[draw=blue, {Stealth}-{Stealth}]
(11, -.6) -- (13, -.6) node[fill=white]{$b+1$} -- (15, -.6);
\end{tikzpicture}
\caption{The parts of the composition $\overline{\varphi(I)}$
for $I\in\mathcal A$, and the value of $v$ in $\overline{\varphi(I)}$.}
\label{fig:reverse.phiI.a}
\end{figure}
Now we are in a position to simplify the expression of $d_I'$.
On one hand, for any $m\ge 0$, 
\begin{align*}
\sum_{k=1}^c \Theta_I^+(k)
&=\brk1{(i_1-1)+\dots+1}
+\dots
+\brk1{(i_m-1)+\dots+1}
+\brk1{(i_{m+1}-1)
+\dots
+(i_{m+1}-u)}
\\
&=
\binom{i_1}{2}
+\dots
+\binom{i_m}{2}
+i_{m+1}u-\binom{u+1}{2}.
\end{align*}
By \cref{pf:d>=0:def:u}, we find
\begin{align*}
\sum_{k=2}^c \Theta_I^+(k)
&=
-
\Theta_I^+(1)
+
\binom{c-u}{2}
-e_2(i_1,\,\dots,\,i_m)
+i_{m+1}u
-\binom{u+1}{2}
\\
&=
-(i_1-1)
+
\binom{c}{2}-(c-i_{m+1})u-e_2(i_1,\,\dots,\,i_m),
\end{align*}
in which $e_2(i_1,\,\dots,\,i_m)=0$ for $m\le 1$.
In the same fashion, 
by \cref{pf:d>=0:def:v}, one may infer that
\[
\sum_{k=a}^{a+c-2}\Theta_{\overline{\varphi(I)}}^-(k)
=\sum_{j=2}^t\binom{i_j}{2}
+\binom{v+1}{2}
=\binom{c-1}{2}-e_2(i_2,\,\dots,\,i_t,\,v+1).
\]
Therefore,
\begin{align*}
d_I'
&=\sum_{k=2}^{c}
\Theta_{I}^{+}(k)
-\sum_{k=a}^{a+c-2}
\Theta_{\overline{\varphi(I)}}^-(k)\\
&=\binom{c}{2}
-(c-i_{m+1})u
-(i_1-1)
-e_2(i_1,\,\dots,\,i_m)
-\binom{c-1}{2}
+e_2(i_2,\,\dots,\,i_t,\,v+1)\\
&=
c-cu+i_{m+1}u-i_1
+e_2(i_2,\,\dots,\,i_t,\,v+1)
-e_2(i_1,\,\dots,\,i_m).
\end{align*}
If $m=0$, then $i_1\ge c$ and
the above expression reduces to
\[
d_I'
=(i_1-c)(u-1)
+e_2(i_2,\,\dots,\,i_t,\,v+1),
\]
which is nonnegative,  since $u\ge 1$ 
and since any evaluation of the $e_2$-function is nonnegative.
Below we can suppose that $m\ge 1$.
Recall that $t\ge m$. 
We evaluate the $e_2$-functions as
\AddToHook{env/align/begin}{\crefalias{equation}{equation}}
\begin{align}
\label{pf:e2.1}
e_2(i_1,\,\dots,\,i_m)
&=
e_2(i_2,\,\dots,\,i_m)
+(i_2+\dots+i_m) i_1,
\quad\text{and}\\
\label{pf:e2.2}
e_2(i_2,\,\dots,\,i_t,\,v+1)
&=
e_2(i_2,\,\dots,\,i_m)
+(i_2+\dots+i_m)
\brk1{i_{m+1}+\dots+i_t+v+1}
+y,
\end{align}
where
$y=e_2(i_{m+1},\,\dots,\,i_t,\,v+1)$.
By \cref{pf:d>=0:def:u,pf:d>=0:def:v}, 
it is direct to see that
\begin{align}
\label{pf:i2im}
i_2+\dots+i_m
&=c-u-i_1
\quad\text{and}\\
\label{pf:v-u}
i_{m+1}+\dots+i_t
&=i_1+u-v-2.
\end{align}
Taking the difference of \cref{pf:e2.1,pf:e2.2} and 
using \cref{pf:i2im,pf:v-u}, we can deduce that
\[
e_2(i_2,\,\dots,\,i_t,\,v+1)
-e_2(i_1,\,\dots,\,i_m)
=
(c-u-i_1)(u-1)
+y.
\]
It follows that 
\[
d_I'
=c-cu+i_{m+1}u-i_1+(c-u-i_1)(u-1)+y
=Cu+y, 
\]
where $C=i_{m+1}+1-i_1-u$.
If $C\ge 0$, then $d_I'\ge 0$, 
since $u,y\ge 0$.

Below we can suppose that $C<0$.
If $m=t$, then \cref{pf:v-u} reduces to $v-u=i_1-2$.
Since $t=m\le p-1$, 
we can derive from \cref{ub:v:t<=p-1} that
$C=i_{m+1}-v-1=i_{t+1}-v-1\ge 0$,
a contradiction. Thus $t\ge m+1$,
which allows us to rewrite $y$ as
\[
y
=
e_2(i_{m+2},\,\dots,\,i_t,\,v+1)
+
i_{m+1}
(i_{m+2}+\dots+i_t+v+1).
\]
By \cref{pf:v-u}, 
\[
i_{m+2}+\dots+i_t
=i_1+u-v-2-i_{m+1}
=-C-v-1.
\]
It follows that 
\[
d_I'
=Cu
+
e_2(i_{m+2},\,\dots,\,i_t,\,v+1)
-
i_{m+1}
\cdotp 
C
=e_2(i_{m+2},\,\dots,\,i_t,\,v+1)
-C(i_{m+1}-u).
\]
By \cref{ub:u}, we deduce that $d_I'\ge 0$.
This completes the proof.
\end{proof}

\section{All clocks are $e$-positive}\label{sec:clock}

%: def:psqt:Delta
In this section, we restrict ourselves to $c=2$.
The clock $G=\theta_{ab2}$ has order $n=a+b+1$.
Let $I=i_1\dotsm i_z$ be a composition of $n$ such that $w_I>0$.
Recall from \cref{def:Delta} that
\begin{equation}\creflabel[def]{def:Delta.b+1}
\Delta_I(b+1)
=\begin{dcases*}
s(i_p-s-i_1),
& if $i_1\le i_p-s$,\\
e_2(i_p-s,\,i_{p+1},\,\dots,\,i_q,\,t),
& otherwise,
\end{dcases*}
\end{equation}
where $p$ and $s$ are defined by the conditional equation
\begin{equation}\creflabel[def]{def:ps}
b+1
=i_1+\dots+i_{p-1}+s
\quad\text{with $1\le p\le z$ and $1\le s\le i_p$}, 
\end{equation}
and $q$ and $t$ are defined by the conditional equation
\begin{equation}\creflabel[def]{def:qt}
b+1
=i_2+\dots+i_q+t
\quad\text{with $1\le q\le z$, $1\le t\le i_{q+1}$ and $i_{z+1}=i_1$}.
\end{equation}
For any $0\le j\le z-p$,
we write
\[
L_j
=
i_1\dotsm i_{p+j}
\quad\text{and}\quad
R_j=i_{p+j+1}\dotsm i_z.
\]
Define $L_I=L_0$ and $R_I=R_0$. 
Then by \cref{lem:psqt},
\begin{equation}\label{Theta-}
\Theta_{\overline{I}}^-(a)
=a-\abs{R_I}
=i_p-s.
\end{equation}
Recall from \cref{def:DI} that
$
D_I
=\Theta_I^+(2)
-\Theta_{\overline{\varphi(I)}}^-(a)
+\Delta_I(b+1)$.

\begin{lemma}\label{lem:D.i1=1}
Let $I\vDash n$ such that $w_I>0$. If $i_1=1$, then $D_I\ge 0$,
where $D_I$ is defined by \cref{def:DI}.
\end{lemma}
\begin{proof}
Let $I=i_1\dotsm i_z$. 
In view of the form of $D_I$,
we will evaluate $\Theta_I^+(2)$, 
$\Theta_{\overline{\varphi(I)}}^-(a)$,
and 
$\Delta_I(b+1)$, respectively.
Since $i_1=1$, we have 
$\Theta_I^+(2)=i_2-1$.

If $I\in\mathcal A$, then $D_I\ge 0$ by \cref{lem:A}.
Alternatively, we have $i_p-s=0$, $\varphi(I)\in\mathcal A$ 
and $\Theta_{\overline{\varphi(I)}}^-(a)=0$.
By \cref{lem:psqt},
$i_1=i_{p+1}+\dots+i_q+t$.
Since $i_1=1$, we find $q=p$. 
By \cref{def:Delta.b+1}, 
one may compute $\Delta_I(b+1)=e_2(0,t)=0$. Hence $D_I=i_2-1\ge 1$.

Suppose that $I\not\in\mathcal A$. 
By \cref{lem:psqt},
\begin{equation}\creflabel[ineq]{range.s}
1\le s\le i_p-1.
\end{equation}
It follows that $p\ge 2$.
By \cref{def:phi}, we find
$\varphi(I)
=i_1 \cdot 
i_{p-1}\dotsm i_2
\cdot i_p R_I$.
Since $\abs{R_I}<a<i_p+\abs{R_I}$,
\begin{equation}\label{pf:Theta.varphi.reverse.a}
\Theta_{\overline{\varphi(I)}}^-(a)
=a-\abs{R_I}
=i_p-s.
\end{equation}
By \cref{def:Delta.b+1}, 
$\Delta_I(b+1)=s(i_p-s-i_1)$.
Therefore,
\[
D_I
=
(i_2-1)
-(i_p-s)
+s(i_p-s-1)
=
(s-1)
(i_p-s-1)
+
(i_2-2), 
\]
which is nonnegative by \cref{range.s}.
\end{proof}

%: def:RQR
\Cref{lem:D.i1=1} leads us to the compositions with all parts at least $2$. Let
\[
\mathcal W=\{I\vDash n\colon i_1,i_2,\dots\ge 2\}.
\]
Now we introduce the aforementioned partial reversal transformation~$\psi$.
For any composition $I\in\mathcal W$, define
\begin{equation}\creflabel[def]{def:psi}
\psi(I)=\overline{L_I}R_I.
\end{equation}
It is clear that $\psi(\mathcal A)\subseteq\mathcal A$.
As will be seen, 
this transformation 
induces a grouping of the compositions in $\mathcal W$
such that the sum of the functions $D_Iw_Ie_I$ is 
$e$-positive in each group. Let
\[
\mathcal W_>
=
\brk[c]1{I\in\mathcal W
\colon
i_1>\Theta_{\overline{I}}^-(a)
}
\quad\text{and}\quad
\mathcal W_\le
=
\brk[c]1{I\in\mathcal W
\colon 
i_1\le \Theta_{\overline{I}}^-(a)
}.
\]
Then $\mathcal W=\mathcal W_>\sqcup \mathcal W_\le$.
By \cref{lem:psqt}, these sets can be defined alternatively as
\[
\mathcal W_>
=
\brk[c]1{I\in\mathcal W
\colon
q\ge p
}
\quad\text{and}\quad
\mathcal W_\le
=
\brk[c]1{I\in\mathcal W
\colon 
q=p-1
}.
\]
We claim that
\[
\psi(\mathcal W_\le)\subseteq \mathcal W_>.
\]
In fact,
let $I\in \mathcal W_\le$ and $J=\psi(I)$.
Then $J=i_p\dotsm i_1 R_I$.
Since $i_1\le \Theta_{\overline{I}}^-(a)$,
we find $\abs{R_J}\ge i_1+\abs{R_I}$.
By \cref{def:ps} of $p$,
we can infer that
\[
i_p
>a-\abs{R_I}
>a-\abs{R_I}-i_1
\ge a-\abs{R_J}.
\]
By \cref{Theta-}, we obtain $J\in \mathcal W_>$.
This proves the claim.
As a consequence, 
\begin{equation}\label{X0123}
X_G
=
\sum_{I\vDash n,\ i_1=1}
D_I
w_I
e_I
+
\sum_{I\in \mathcal W_>}
c_I
e_I,
\end{equation}
where 
\[
c_I
=
D_I
w_I
+
\sum_{H\in\psi^{-1}(I)\cap \mathcal W_\le}
D_H
w_H.
\]

Next, we are going to figure out the fiber $\psi^{-1}(I)$
(see \cref{prop:psi.inverse.i1<=})
and to estimate $D_I$ (see \cref{lem:D.W>,lem:D.W<=,lem:DUr<0}).
We present an interpretation of the difference $q-p$,
which will be of use in the proofs of \cref{prop:psi.inverse.i1<=,lem:D.W>}.

\begin{proposition}\label{prop:q-p}
If $I=i_1\dotsm i_z\in \mathcal W_>$, then 
$
q-p=\max\{
j\le z-p\colon
\abs{R_j}>a-i_1
\}$.
\end{proposition}
\begin{proof}
Let $I\in \mathcal W_>$. 
Let $m$ be the desired maximum. 
Then $m\ge 0$ since $i_1>\Theta_{\overline{I}}^-(a)=a-\abs{R_0}$.
If $i_1>a$, then $m=z-p$ as defined.
By \cref{lem:psqt}, we find $q=z$ and $m=q-p$ as desired.
Suppose that $i_1\le a$. 
By \cref{lem:psqt}, we find $q\le z-1$.
Then 
both $R_{q-p}$ and $R_{q-p+1}$
are well defined, and
\[
\abs{R_{q-p}}
=i_{q+1}+\dots+i_z
=a-i_1+t
>a-i_1
\ge a-i_1+t-i_{q+1}
=
\abs{R_{q-p+1}}.
\]
This proves $m=q-p$ as desired.
\end{proof}

% def:U
For any $I\in\mathcal W$ and any $0\le r\le q-p$, define $H_r=\overline{L_r}R_r$.
Then $e_{H_r}=e_I$.

\begin{proposition}\label{prop:psi.inverse.i1<=}
If $I\in \mathcal W_>$,
then $R_{H_r}=R_r$ for all $1\le r\le q-p$, and 
\[
\psi^{-1}(I)\cap \mathcal W_\le
=\{H_1,\,\dots,\,H_{q-p}\}.
\]
\end{proposition}
\begin{proof}
Let $I\in \mathcal W_>$ and $J\in \psi^{-1}(I)\cap \mathcal W_\le$.  
Then
\[
\overline{L_J}R_J
=\psi(J)
=I
=i_1\dotsm i_p R_I.
\]
Since
$\abs{R_J}\le a$
and 
$\abs{R_I}\le a<i_p+\abs{R_I}$,
we deduce that $R_J$ is a suffix of $R_I$,
namely 
$R_J=R_r$ for some $0\le r\le z-p$.
It follows that
\[
J=\overline{L_r}R_r=H_r=i_{p+r}\dotsm i_1R_r.
\]
Since $R_r=R_J$ is the longest suffix of $J$ whose size is at most~$a$, 
we find $i_1+\abs{R_r}>a$. By \cref{prop:q-p}, we find $r\le q-p$.
On the other hand, if $r=0$, then 
$j_1=i_p>\Theta_{\overline{I}}^-(a)=\Theta_{\overline{J}}^-(a)$,
contradicting $J\in \mathcal W_\le$. Thus $r\ge 1$.
This proves $\psi^{-1}(I)\cap \mathcal W_\le
\subseteq\{H_1,\,\dots,\,H_{q-p}\}$.

Conversely, let $1\le r\le q-p$ and 
$J=H_r=i_{p+r}\dotsm i_1 R_r$.
By \cref{prop:q-p},
\[
\abs{R_r}
<\abs{R_0}
\le a
<i_1+\abs{R_r}.
\]
Thus $R_J=R_r$ and
$\psi(J)
=L_r R_r
=I$.
Since 
\[
j_1
=i_{p+r}
=\abs{R_{r-1}}-\abs{R_r}
\le \abs{R_0}-\abs{R_r}
\le a-\abs{R_r}
=\Theta_{\overline{J}}^-(a),
\]
we find $J\in \mathcal W_\le$.
This proves $\{H_1,\,\dots,\,H_{q-p}\}
\subseteq \psi^{-1}(I)\cap \mathcal W_\le$, and completes the whole proof.
\end{proof}

Now we estimate the values of $D_I$.
We proceed for $I$ in the sets $\mathcal W_>$ and $\mathcal W_\le$, respectively.

%: lem:D.W>
\begin{lemma}\label{lem:D.W>}
Let $I\in\mathcal W_>$. 
Then $D_I\ge i_1-2\ge 0$. 
If $q\ge p+1$ in addition,
then we have the following.
\begin{enumerate}
\item
If $I\in\mathcal A$,
then $i_1\ge 3$ and $D_I\ge 2i_1-3$.
\item
If $I\in\mathcal W_>\backslash \mathcal A$,
then $i_1\ge 4$ and $D_I\ge i_1+2$. 
\end{enumerate}
\end{lemma}
\begin{proof}
Let $I=i_1\dotsm i_z\in\mathcal W_>$ and $s'=i_p-s$.
Since $i_1\ge 2$ and $i_1>s'$,
we have 
\[
\Theta_I^+(2)=i_1-2
\quad\text{and}\quad
\Delta_I(b+1)=e_2(s',\,i_{p+1},\,\dotsm,\,i_q,\,t).
\]
As in the proof of \cref{lem:D.i1=1},
one may show \cref{pf:Theta.varphi.reverse.a}. 
Then by \cref{def:DI} of $D_I$, 
\begin{equation}\label{pf:D=.W>}
D_I
=i_1-2-s'+e_2(s',\,i_{p+1},\,\dotsm,\,i_q,\,t).
\end{equation}
Since $t\ge 1$, we find $e_2(s',\,i_{p+1},\,\dotsm,\,i_q,\,t)\ge s'$,
and $D_I\ge i_1-2$ as desired.

Below we suppose that $q\ge p+1$.
Since $i_{p+1}\ge 2$, 
by \cref{lem:psqt}, 
\begin{equation}\creflabel[ineq]{pf:2<=i1-s'-t}
2
\le 
i_{p+1}+\dots+i_q
=i_1-s'-t.
\end{equation}
Since $s'\ge 0$ and $t\ge 1$,
we can infer that $i_1\ge 3$, $1\le t\le i_1-2$, and
\AddToHook{env/align/begin}{\crefalias{equation}{ineq}}
\begin{align}
\notag
e_2(s',\,i_{p+1},\,\dotsm,\,i_q,\,t)
&=
e_2(i_{p+1},\,\dots,\,i_q)
+(s'+t)(i_{p+1}+\dots+i_q)
+s't
\\
\creflabel[ineq]{pf:e2.W>}
&\ge 
(s'+t)(i_1-s'-t)
+s't.
\end{align}
Now we are ready to give a lower bound for $D_I$.

\begin{enumerate}
\item
If $I\in\mathcal A$, 
then $s'=0$.
Note that $i_{p+1},\dots,i_q\ge 2$ and $1\le t\le i_1-2$.
By \cref{pf:D=.W>,pf:e2.W>},
\[
D_I
\ge i_1-2+t(i_1-t)
\ge 2i_1-3.
\]
In fact, 
since the positive numbers $t$ and $i_1-t$ have constant sum $i_1$,  
their product attains the minimum when their difference attains the maximum, i.e., when $t=1$.
\item
If $I\not\in\mathcal A$, then $s'\ge 1$.
Since $t\ge 1$, we can infer 
by \cref{pf:D=.W>,pf:e2.W>,pf:2<=i1-s'-t} that
\[
D_I
\ge i_1-2-s'
+2(s'+t)
+s't
\ge i_1-2-s'
+2(s'+1)
+s'
\ge i_1+2,
\]
as desired.
On the other hand, by \cref{Theta-,prop:q-p},
\[
a
\ge \abs{R_I}
=i_{p+1}+\abs{R_1}
\ge 2+(a-i_1+1).
\]
If $i_1=3$, then the inequality above implies $\abs{R_I}=a$, contradicting $I\not\in\mathcal A$. Thus $i_1\ge 4$.
\end{enumerate}
This completes the proof.
\end{proof}

The lower bounds for $D_I$ are sharp in sense of the following examples.
\begin{enumerate}
\item
For $(a,b)=(2,2)$ and $I=32$,
we have $I\in\mathcal A$, $H_1=23$ and $D_I=3$.
\item
For $(a,b)=(3,2)$ and $I=42$,
we have $I\in\mathcal W\backslash \mathcal A$, $H_1=24$ and $D_I=6$.
\end{enumerate}

We then estimate the number $D_I$ for $I\in\mathcal W_\le$.
%: lem:D.W<=
\begin{lemma}\label{lem:D.W<=}
Let $I\in\mathcal W_\le$.
Then 
\[
D_I=(s-1)(i_p-s-i_1)-2\ge -2.
\]
Moreover, if $D_I<0$, then $s\in\{1,\,2,\,i_p-i_1\}$.
\end{lemma}
\begin{proof}
Let $I\in \mathcal W_\le$.
Since $i_1\ge 2$, we have 
$\Theta_I^+(2)=i_1-2$.
Note that \cref{pf:Theta.varphi.reverse.a} still holds, even if $p=1$.
Since $i_1\le i_p-s$, \cref{def:Delta.b+1} gives $\Delta_I(b+1)=s(i_p-s-i_1)$.
Therefore,
\[
D_I
=
(i_1-2)
-(i_p-s)
+s(i_p-s-i_1),
\]
which coincides with the desired one.
It is at least $-2$ since $1\le s\le i_p-i_1$.
When $D_I<0$, 
it is easy to derive the desired range of $s$ by 
letting $D_I=-2$ and $D_I=-1$ respectively.
\end{proof}

In order to show $c_I\ge 0$,
we track the preimages of $H_r$ that have negative $D$-values.

\begin{lemma}\label{lem:DUr<0}
Let $I\in \mathcal W_>$.
Suppose that $q-p\ge 1$.
If $D_{H_r}<0$ for some $1\le r\le q-p$,
then either (i) 
$r=q-p$, or (ii)
$I\in\mathcal A$ and $r=1$.
\end{lemma}
\begin{proof}
Suppose that $D_J<0$, 
where 
$J=H_r=i_{p+r}\dotsm i_1 R_r$
for some $1\le r\le q-p$.
By \cref{prop:psi.inverse.i1<=}, we have $R_J=R_r$.
Let $s_r$ be the number $s$ defined by \cref{def:ps},
in which the composition $I$ is replaced with $H_r$.
Applying \cref{Theta-} to $J$, we obtain
\begin{equation}\label{pf:sJ}
s_r=i_1-a+\abs{R_r}.
\end{equation}
Since $\abs{R_J}=\abs{R_r}\le\abs{R_1}<\abs{R_0}\le a$,
we find $J\not\in\mathcal A$.
By \cref{lem:D.W<=}, we find $s_r\in\{1,\,2,\,i_1-i_{p+r}\}$.

If $s_r=i_1-i_{p+r}$, 
then \cref{pf:sJ} implies
$a=i_{p+r}+\abs{R_r}=\abs{R_{r-1}}\le\abs{R_0}\le a$.
Since the equality holds, we find $r=1$ and $\abs{R_0}=a$,
i.e., $J=H_1$ and $I\in\mathcal A$.
Now suppose that $s_r\in\{1,2\}$. Assume $r<q-p$.
Then $\abs{R_r}\ge \abs{R_{r+1}}+2\ge \abs{R_{q-p}}+2$.
By \cref{pf:sJ}, we deduce that $s_r\ge s_{q-p}+2\ge 3$, a contradiction. 
This proves $r=q-p$, and completes the whole proof.
\end{proof}

For example, for $(a,b)=(6,4)$, we have the following.
\begin{itemize}
\item
For $I=722$, 
we have 
$I\in \mathcal W_>$, 
$q-p=2$,
$H_1=272$,
$H_2=227$, 
$D_{H_1}=2$
and 
$D_{H_2}=-2$.
\item
For $I=5222$, 
we have 
$I\in \mathcal W_>\cap\mathcal A$,
$q-p=2$,
$H_1=2522$, 
$H_2=2252$,
and
$D_{H_1}=D_{H_2}=-2$.
\end{itemize}

Now we are in a position to prove \cref{thm:X.clock}.
\begin{proof}[Proof of \cref{thm:X.clock}]
By \cref{X0123,lem:D.i1=1}, it suffices to show that $c_I\ge 0$ for each $I\in\mathcal W_>$.
Let $I\in\mathcal W_>$. Then $q\ge p$ by \cref{prop:q-p}.
If $q=p$, then the sum vanishes by \cref{prop:psi.inverse.i1<=}.
In this case,  $c_I=D_Iw_I\ge 0$ by \cref{lem:D.W>}.
Below we suppose that $q\ge p+1$.
\begin{enumerate}
\item
If $I\not\in\mathcal A$,
then by \cref{lem:D.W>,lem:D.W<=,lem:DUr<0}, 
we have $i_1\ge 4$ and 
\[
c_I
\ge 
D_I
w_I
+
D_{H_{q-p}}
w_{H_{q-p}}
\ge 
w_{1I}
\brk3{
\frac{i_1(i_1+2)}{i_1-1}
-\frac{2i_q}{i_q-1}
}
\ge 
2w_{1I}
\brk3{
\frac{3}{i_1-1}
+2
-\frac{1}{i_q-1}
}
>0.
\]
\item
If $I\in\mathcal A$, 
then $i_1\ge 3$ and $D_I\ge 2i_1-3$ by \cref{lem:D.W>}.
If $q-p=1$, then 
by \cref{lem:D.W<=,lem:DUr<0},
\begin{align*}
c_I
\ge 
D_I
w_I
+
D_{H_1}
w_{H_1}
&\ge 
w_{1I}
\brk3{
\frac{i_1(2i_1-3)}{i_1-1}
-\frac{2i_{p+1}}{i_{p+1}-1}
}\\
&=
w_{1I}
\brk3{
2i_1-3-\frac{1}{i_1-1} 
-\frac{2}{i_{p+1}-1}
}
\ge 
w_{1I}
\brk2{
3-\frac{1}{2} 
-2
}
>0.
\end{align*}
Otherwise $q-p\ge 2$. In this case,
\cref{prop:q-p} implies
\[
a=\abs{R_I}=i_{p+1}+i_{p+2}+\abs{R_2}
\ge 2+2+(a-i_1+1).
\]
Thus $i_1\ge 5$.
By \cref{lem:D.W<=,lem:DUr<0},
\begin{align*}
c_I
\ge 
D_I
w_I
+
D_{H_1}
w_{H_1}
+
D_{H_{q-p}}
w_{H_{q-p}}
&\ge 
w_{1I}
\brk3{
\frac{(2i_1-3)i_1}{i_1-1}
-\frac{2i_{p+1}}{i_{p+1}-1}
-\frac{2i_{q}}{i_{q}-1}
}
\\
&=
w_{1I}
\brk3{
2i_1-5-\frac{1}{i_1-1} 
-\frac{2}{i_{p+1}-1}
-\frac{2}{i_{q}-1}
}\\
&\ge 
w_{1I}
\brk3{
5-\frac{1}{4} 
-2
-2
}
>0.
\end{align*}
\end{enumerate}
This completes the proof.
\end{proof}

\bibliographystyle{abbrvnat}
\bibliography{csf.bib}

\end{document}